\documentclass[12pt]{amsart}
\usepackage{amssymb,amsmath,txfonts,amsrefs, amscd}
\usepackage{array,multirow,bigstrut}
\usepackage{appendix,mathtools}
\usepackage[all]{xy}
\usepackage{xargs}  
\usepackage[pdftex,dvipsnames]{xcolor}  
\usepackage[colorinlistoftodos,prependcaption,textsize=tiny]{todonotes}
\usepackage{hyperref} 

\newcommandx{\unsure}[2][1=]{\todo[linecolor=red,backgroundcolor=red!25,bordercolor=red,#1]{#2}}

\newcommandx{\change}[2][1=]{\todo[linecolor=blue,backgroundcolor=blue!25,bordercolor=blue,#1]{#2}}

\newcommandx{\info}[2][1=]{\todo[linecolor=OliveGreen,backgroundcolor=OliveGreen!25,bordercolor=OliveGreen,#1]{#2}}

\newcommandx{\improvement}[2][1=]{\todo[linecolor=Plum,backgroundcolor=Plum!25,bordercolor=Plum,#1]{#2}}
\newcommandx{\thiswillnotshow}[2][1=]{\todo[disable,#1]{#2}}

\newcommand{\closure}[2][3]{%
{}\mkern#1mu\overline{\mkern-#1mu#2}}

\newtheorem{thm}{Theorem}[section]
\newtheorem{prop}[thm]{Proposition}
\newtheorem{lemma}[thm]{Lemma}
\newtheorem{remark}[thm]{Remark}

\newtheorem{defi}[thm]{Definition}

\DeclareMathOperator{\Ima}{Im}
\numberwithin{equation}{section}
     
\def\S{Steklov\,}

\def\HM{W^{k-\frac{1}{p},\,p}(\partial M^n)}
\def\hm{W^{k-1-\frac{1}{p},\,p}(\partial M^n)}
\def\wm{W^{k,\,p}(M^n)}

\def\DN{the Dirichlet-to-Neumann operator\,}
\def\R{\mathbb{R}}

\def\L{\Lambda}

\def\Gk{\mathcal{G}^{k}}

\def\Ck{\mathcal{C}^{k}(\closure{M^n})}

\def\Bn{\mathbb{B}^n}

\begin{document}

\title[]{Generic properties of Steklov eigenfunctions}

\author{Lihan Wang}
\address{Department of Mathematics and Statistics, CSULB }\email{lihan.wang@csulb.edu}

\date{\today}

\keywords{}

\begin{abstract}

Let $M^n$ be a smooth compact manifolds with smooth boundary. We show that for a generic $C^k$ metic on $\closure{M^n}$ with $k>n-1$, the nonzero Steklov eigenvalues are simple. Moreover, we also prove that the non-constant Steklov eigenfunctions have zero as a regular value and are Morse functions on the boundary for such generic metric. These results generalize the celebrated results on Laplacians by Uhlenbeck to the Steklov setting.

\smallskip

\end{abstract}

\maketitle
\pdfbookmark[0]{}{beg}

\section{Introduction}
Given a smooth compact Riemannian manifold $M^n$ with smooth boundary $\partial M^n$, \DN, denoted by $\Lambda$, takes a function $f$ on $\partial M^n$ to the normal derivative of the harmonic extension of $f$.
Its eigenvalues and eigenfunctions are often called the \S eigenvalues and Steklov eigenfunctions respectively. Since they were first discussed by Steklov (\cite{S}) in 1902 motivated by physics: these eigenfunctions represent the steady state temperature on bounded domains such that the flux on the boundary is proportional to the temperature. 

There have been intense study of Steklov eigenvalues and eigenfunctions. 
Recently, Fraser and Schoen made great progress in the extremal Steklov eigenvalue problems for general surfaces, higher dimensional domains and the higher Steklov eigenvalues (\cite{FS11}, \cite{FS13}, \cite{FS16}, \cite{FS19}, \cite{FS20}). In particular, in the seminal work \cite{FS16}, Fraser and Schoen revealed a deep connection between the extremal Steklov eigenvalue problems and the free boundary minimal surface theory in the unit Euclidean ball $\mathbb{B}^n$.
 A free boundary minimal surface in $\Bn$ is a surface contained in the unit ball with zero mean curvature which meets the boundary of the ball orthogonally. Such a surface arises variationally as a critical point of the area among surfaces in $\Bn$ with boundaries lying on $\partial \Bn$ but free to vary on $\partial \Bn$. 
Fraser and Schoen showed that if a metric $g$ realizes the maximum of the first non-zero normalized Steklov eigenvalue, then this metric is $\theta$-homothetic to the induced metric on a free boundary minimal surface in $\Bn$. In fact, there exists a branched conformal immersion into $\Bn$ by first Steklov eigenfunctions whose image is a free boundary minimal surface. 

 It would be very helpful if we can understand the Steklov eigenvalues and eigenfunctions better in the study of the free boundary minimal surfaces. And it could have applications in other areas, like  the inverse conductivity problems (\cite{C}, \cite{SU}) and cloaking (\cite{GKLU}) where \DN has been used. See \cite{GP} for more geometric questions related to Steklov eigenvalues and \cite{GKLP} for recent results about Steklov eigenvalues.

Our purpose of this paper is to show certain generic properties of Steklov eigenfunctions. In her seminal work \cite{U}, Uhlenbeck established generic properties of eigenfunctions of second-order elliptic differential operators on compact manifolds. In particular, one of celebrated results asserts that  nonzero eigenvalues of the Laplacian $\Delta_g$ are simple for generic metrics. Many works have been done to generalize this result in different directions. For instance, the nonzero eigenvalues of the Hodge Laplacian are shown to be generically simple on a closed manifold of dimension 3 in \cite{EP12}. In addition to second-order differential operators, eigenvalues of the Dirac operator on a spin three manifold are also shown to be simple for generic metrics in \cite{DM}. The case of conformally covariant operators was discussed in \cite{CY}.

Lastly, we note that the spirit of Uhlenbeck's work has been influential in other non-linear geometric variational problems. In particular, in the recent celebrated developments of the min-max theory of minimal hypersurfaces and related volume spectrum, the Multiplicity One Conjecture (\cite{ZX}) can be viewed an analogous result to Uhlenbeck's work. The free boundary version is treated in \cite{SWZ}.


In this paper, we prove that nonzero \S eigenvalues are simple for generic metrics. Moreover, we prove the corresponding Steklov eigenfunctions generically have zero as a regular value and are Morse functions on the boundary. That is, the nodal sets and the sets of critical points of Steklov eigenfunctions on the boundary consist of isolated points for generic metrics. The study of nodal sets and critical points of eigenfunctions
is one of oldest topic in geometric spectral theory. In the case of Steklov eigenfunctions, the nodal sets and critical points remain largely unexplored. It would be of interest to see what further results these generic properties may lead to. After the completion of the work, we were informed that the generic simplicity of \S eigenvalues was conjectured in \cite{GP}.


Our main result can be stated as follows:
\vspace{- 2.5mm}
\begin{thm}[Main Theorem]
Let $M^n$ be a smooth compact manifold with smooth boundary $\partial M^n$ for $n\geq 2$. Given an integer $k >n-1$, there exists a residual set of Riemannian metrics which are $C^k$ on $\closure{M^n}$, for which 
 the following properties hold at non-constant \S eigenfunctions:
\vspace*{-1.5mm}
\begin{enumerate}
\item The eigenspaces of nonzero \S eigenvalues are one dimensional, i.e., nonzero \S eigenvalues are simple.
\item Zero is a regular value of eigenfunctions restricted to $\partial M^n$.
\item The eigenfunctions are Morse functions on $\partial M^n$.
\end{enumerate}
 \end{thm}

 Following the spirit of Uhlenbeck's work in \cite{U}, we employ the infinite dimensional transversality theory to demonstrate these generic properties for \S eigenfunctions. As known, \DN is a first-order, elliptic, self-adjoint and pseudo-differential operator. Various essential difficulties arise in the extension of Uhlenbeck's proofs and results. One challenge appears in the verification of regular values of evaluation maps. In the case of differential operators, Green's function is employed in \cite{U}. However this approach does not work in the case of \DN. Since \DN has a compact resolvent, we take a different approach using resolvent operators for this purpose. Technical details are listed in Section 3. It is worth mentioning that resolvent operators are also used to prove generic properties of eigenvalues and eigenfunctions of Hodge Laplacians in \cite{EP12}.


 
Due to the non-local property of \DN, some technical complications arise in deriving the density lemma (Lemma \ref{density}) which plays a key role in the application of transversality theorems. In the case of differential operators, the density lemma can be derived directly by variation formulas of operators with respect to the underlying metric, which are local. However, in our case, the variation formula of $\L f$ with respect to the metric depends on not only the variation of $g$ on $\partial M^n$, but also the harmonic extension of $f$ and the corresponding variation on $M^n$. So the usual approach does not work for \DN. To overcome the complication, we use variation formulas both of $\L f$ and the the harmonic extension of $f$, and apply Green's formula carefully to derive the desired density lemma (Lemma \ref{density}). Moreover, this density lemma provides the existence of certain conformal variations of the metric which are needed in proofs in Section 3. The precise discussion is included in Section 2 and Section 3.

The article is organized as follows. In Section 2, we set up the definition of $\Lambda$ and evaluate its variation with respect to the metric. We also derive the variation of harmonic extensions with respect to the metric. The weak continuation principle of Steklov eigenfunctions is established in this section too. In Section 3, we recall basic definitions and transversality theorems. Then we prove three generic properties of Steklov eigenfunctions in Main Theorem separately through Theorem \ref{Gen1} (Section 3.1), Theorem \ref{Gen2}(Section 3.2) and Theorem \ref{Gen3} (Section 3.2). 
 
The author would like to thank Xin Zhou for suggesting this topic and helpful comments. The author thanks the referee for helpful comments to improve the exposition and pointing out the conjecture of generic simplicity  of \S eigenvalues in \cite{GP}.

\section{Preliminary Lemmas}

Throughout this paper, $M^n$ will stand for a smooth compact manifold with smooth boundary $\partial M^n$ and $n\geq 2$. We shall consider the space $\Gk$ of Riemannian metrics $g$ which are $C^{k}$ on $\closure{M^n}$. For $g \in \Gk$ with integers $k \geq 1$ and $p \geq 2$, \DN is defined as: \begin{align}\label{eq1}
\begin{aligned}
\Lambda: \HM &\rightarrow \hm\\
f &\rightarrow u_n
\end{aligned}
\end{align} with $u$ as the solution of
\begin{align}\label{eq1.1}
\left\{ 
\begin{aligned}
\Delta u&=0,\, && M^n \\ 
u&=f,\, &&\partial M^n.
\end{aligned}%
\right.
\end{align} 
In fact, according to trace theorems (Chapter 2 in \cite{NJ}) and elliptic theory (Chapter 9 in \cite{GT}), for any $f \in \HM$, there exists the unique harmonic extension $u \in \wm$ of $f$ when $g \in \Gk$. And its normal derivative $u_n$ belongs to $\hm$.

\begin{remark} We can make assumptions for this definition weaker. Actually $\Lambda$ above is still well-defined when $\closure{M^n}$ is $C^{k,1}$ and $g$ is  $C^{k-1,1}$ on $\closure{M^n}$ according to results in \cite{NJ} and \cite{GT}.
\end{remark}


In this section, we will evaluate the variations of the harmonic extension $u$ and the \DN $\Lambda$ with respect to the metric $g$. We shall denote the tangent space of $\Gk$ by $T\Gk$, which can be identified as the space of symmetric tensor fields of class $C^k$ on $\closure{M^n}$ and of type $(0,2)$. The variation of the metric $g$ shall be denoted by $h=Dg\in T\Gk$. When the variation is conformal, we denote it as $h=\sigma g$ for some positive function $\sigma \in \Ck$. The symbol $h^{-1}=Dg^{-1}$ refers to the variation of $g^{-1}$. The symbol $D_g$ refers to the differential with respect to $g$.

We use $g_{ij}, g^{ij}, h_{ij}$ and $h^{ij}$ to denote the components of $g, g^{-1}, h, h^{-1}$ in local coordinates respectively. We employ the Einstein summation convention throughout this paper. The trace of a tensor $S$ with respect to $g$ will be written as $\rm{tr}_g S\vcentcolon = g^{ij}S_{ij}$. The volume element of $M^n$ shall be denoted by $dV$ and the area element of $\partial M$ shall be denoted by $dA$.

At first, we will derive the variations of harmonic extensions with respect to the metric $g$.

\begin{lemma}[Variations of harmonic extensions]\label{VH}
Consider a function $f\in \HM$ and its harmonic extension $u \in \wm$ with $k\geq 2$ and $p\geq 2$. Let $h=Dg$ be the variation of $g$. With $f$ fixed, we denote the variation of $u$ with respect to $g$ by $v=D_g u$. Then $v$ satisfies  
\begin{align}\label{eq2}
 \left\{
\begin{aligned}
\Delta v+D_g(\Delta)u&=0, \,&&M^n\\
v&=0,\, &&\partial M^n
\end{aligned}%
\right.
\end{align} with $D_g(\Delta)u=h^{ij}u_{ij}+\partial_{i}(h^{ij})u_j+\frac{1}{2} \rm{tr}(h^{-1}\partial_{i}g+g^{-1}\partial_{i}h)g^{ij}u_j+\frac{1}{2} \rm{tr}(g^{-1}\partial_{i}g)h^{ij}u_j$.

When $h=\sigma g$ is a conformal variation for some $\sigma \in \Ck$, it follows that 
 \begin{align}\label{eq3}
 \left\{
\begin{aligned}
\Delta v-(1-\frac{n}{2})\langle \nabla \sigma, \nabla u\rangle&=0, \, &&M^n\\
v&=0,\, &&\partial  M^n.
\end{aligned}
\right.
\end{align} 
\end{lemma}

\begin{proof}
As the harmonic extension of $f$, $u$ satisfies that
\begin{align}\label{eq4}
 \left\{
\begin{aligned}
\Delta u&=0, \, &&M^n,\\
u&=f,\, &&\partial  M^n.
\end{aligned}
\right.
\end{align}
Do the differential with respect to $g$ on both sides of \eqref{eq4}. It follows that $v=D_g u$ satisfies 
 \begin{align*}
 \left\{
\begin{aligned}
\Delta v+D_g(\Delta)u&=0, \, &&M^n\\
v&=0,\, &&\partial  M^n.
\end{aligned}
\right.
\end{align*}

To calculate $D_g(\Delta$), we express $\Delta$ in terms of the local coordinates $\{x_i\}$:
 \begin{align*} 
\Delta u&= |g|^{-\frac{1}{2}}\partial_{i}\left(|g|^{\frac{1}{2}} g^{ij}u_j\right)=g^{ij}u_{ij}+\partial_{i}(g^{ij})u_j+\frac{1}{2} (\rm{tr}_g\,\partial_{i}g)g^{ij}u_j.
\end{align*} With $u$ fixed, taking the differential with respect to $g$ on both sides yields that
 \begin{align*} 
D_g(\Delta)u&=h^{ij}u_{ij}+\partial_{i}(h^{ij})u_j+\frac{1}{2} (\rm{tr}_h \partial_{i}g+\rm{tr}_g\partial_{i}h)g^{ij}u_j+\frac{1}{2} (\rm{tr}_g \partial_{i}g)h^{ij}u_j.
\end{align*} 

Do the differential on the equality $gg^{-1}=I$. It follows that
\begin{align*}
Dg^{-1}&=-g^{-1}(Dg)g^{-1}\\
&=-g^{-1}hg^{-1}.
\end{align*}  When $h=\sigma g$ is a conformal variation for $\sigma \in \Ck$, it follows that 
\[Dg^{-1}=-\sigma g^{-1}.\]
Plug this into the the above formula of $D_g(\Delta)$, we get 
 \begin{align}\label{eq12}
D_g(\Delta)u&=-\left( \sigma \Delta u+ (1-\frac{n}{2})\langle \nabla \sigma, \nabla u\rangle \right).
\end{align}
Then \eqref{eq3} follows by plugging \eqref{eq4} and \eqref{eq12} into \eqref{eq2}.
\end{proof}

Next, we will derive the variation of $\Lambda$ with respect to the metric $g$ based on the lemma above.

\begin{lemma}[The variation of $\Lambda$]\label{VL}
Consider a function $f\in \HM$ with $k\geq 2$ and $p\geq 2$. Let $u$ denotes its harmonic extension in $\wm$. Let $h=Dg$ be the variation of $g$ and $v=D_g u$ be the variation of $u$ with respect to $g$. Let $(D_g \Lambda)(h)f$ denote the variation of $\Lambda f$ with respect to $g$ in the direction of $h$. Then it satisfies

\begin{align}\label{eq5}
(D_g \Lambda)(h)f&=v_n-h(\nabla u, \vec{n})+\frac{u_n}{2}h(\vec{n}, \vec{n}).
\end{align}Here $\vec{n}$ is the outward norm along $\partial M^n$ and $v_n$ is the derivative of $v$ with respect to $\vec{n}$.

When $h=\sigma g$ is a conformal variation for some $\sigma \in \Ck$, we have 
\begin{align}\label{eq6}
(D_g \Lambda)(h) f&=v_n-\frac{\sigma}{2} u_n.
\end{align}

\end{lemma}

\begin{proof} Fix a point $ x\in \partial  M^n$. In a neighborhood of $x$, there exists some boundary defining function $\rho$ such that $\vec{n}=\frac{\nabla \rho}{|\nabla \rho|}$. 

By the definition \eqref{eq1}, there is
\[\Lambda f= u_{n}=\frac{u_i\rho_j g^{ij}}{|\nabla \rho|}\] with respect to the local coordinates $\{x_i\}$. Taking the differential with respect to $g$ on both sides yields that 
\begin{equation}
\begin{aligned}
(D_g \Lambda)(h) f&=\frac{v_i\rho_j g^{ij}}{|\nabla \rho|}+\frac{u_i\rho_j h^{ij}}{|\nabla \rho|}-\frac{u_n(h^{ij}\rho_i\rho_j)}{2|\nabla \rho|^2}\\
&=v_n+\frac{u_i\rho_j h^{ij}}{|\nabla \rho|}-\frac{u_n(h^{ij}\rho_i\rho_j)}{2|\nabla \rho|^2}.\label{eq2.7}
\end{aligned} 
\end{equation}
Since
\begin{align*}
h(\nabla u, \vec{n})&=-\frac{h^{ij}u_i\rho_j}{|\nabla \rho|}\\
h(\vec{n}, \vec{n})&=-\frac{ h^{ij}\rho_i\rho_j}{|\nabla \rho|^2},
\end{align*} we get \eqref{eq5}.
 
 When $h=\sigma g$, plug $h^{ij}=-\sigma g^{ij}$ into \eqref{eq2.7} :
 \begin{equation*}
 \begin{aligned}
(D_g \Lambda)(h) f&=v_n+\frac{u_i\rho_j h^{ij}}{|\nabla \rho|}-\frac{u_n(h^{ij}\rho_i\rho_j)}{2|\nabla \rho|^2}\\
&=v_n-\sigma\frac{u_i\rho_j  g^{ij}}{|\nabla \rho|}+\sigma\frac{u_n g^{ij}\rho_i\rho_j}{2|\nabla \rho|^2}\\
&=v_n-\sigma u_n+\frac{1}{2}\sigma u_n|\vec{n}|^2\\
&=v_n-\frac{1}{2}\sigma u_n.
\end{aligned} 
\end{equation*} Then we get \eqref{eq6}.
 \end{proof}
\begin{remark}
We can do the variation of $\Lambda$ since it is $C^1$with respect to the metric $g$. In fact, it can be shown that the right hand side of \eqref{eq5} is the directional derivative of $\Lambda$ with respect to $g$ in the direction $h$ using the definition of directional derivative directly. And the the right hand side of \eqref{eq5} depends on $g$ continuously. Thus $\Lambda$ is $C^1$ with respect to $g$.
\end{remark}

In the end, we will prove the weak unique continuation principle for \S eigenfunctions. 

\begin{thm}[Weak unique continuation principle]\label{wucp}
Let $M^n$ be a smooth compact manifold with smooth boundary $\partial M^n$ for $n\geq 2$ and  Riemannian metric $g$ which is $C^k$ on $\closure{M}$. Consider a non-constant function $f\in \HM$ with $k\geq 2$ and $p\geq 2$. Assume that $\Lambda f=\lambda f$. If $f=0$ on a open set of $\partial M^n$, then $f$ vanish on $\partial M^n$.
\end{thm}
\begin{proof}
Let $u$ be the harmonic extension of $f$. Then 
\begin{align}\label{ap1}
\left\{ 
\begin{aligned}
\Delta u&=0,\, && M^n \\ 
u&=f,\, &&\partial M^n\\
u_n&=\lambda f,\, &&\partial M^n.
\end{aligned}%
\right.
\end{align} 
Let $\Gamma$ be an open subset of $\partial M$ where $f=0$. Fix an arbitrary point $x\in \Gamma$. Choose a small neighborhood $\Omega$ of $x$ such that $\Omega\cap\partial M\subset \Gamma$. And there exists a homeomorphism to the unit half ball $B^{+}$ in $\R^n$:
\[
\Phi: \Omega \rightarrow B^{+}=\{(x_1, \cdots, x_n)\in\R^n|\sum x_i^2 \leq 1, x_1\geq 0\}
\] which maps the $x$ to the origin and each boundary point in $\Omega \cap \partial M$ to a point with $x_1=0$. 

With respect to the local coordinates given by $\Phi$, \eqref{ap1} together with the assumption $f=0$ on $\Gamma$ imply that:
\begin{align}\label{ap2}
\left\{ 
\begin{aligned}
|g|^{-\frac{1}{2}}\partial_{i}\left(|g|^{\frac{1}{2}} g^{ij}u_j\right)&=0,\, && B^{+}\\ 
u=f&=0,\, &&\{x_1=0\} \cap B^{+} \\
\frac{\sum g^{ij}\rho_i u_j}{\sqrt{\sum g^{ij}\rho_i\rho_j}}=\lambda f&=0,\, &&\{x_1=0\} \cap B^{+}.
\end{aligned}
\right.
\end{align} 
Here $\rho$ is the boundary defining function around $x$ with $\rho|_{\partial M}=0, |\nabla \rho|_{\partial M}\neq 0$.

Since $u=0$ and $\rho=0$ on the boundary $\{x_1=0\} \cap B^{+}$, we get $u_i=0, \rho_i=0$ for $i=2, \cdots, n$ on $\{x_1=0\} \cap B^{+}$. At the same time, $|\nabla \rho|^2=\sum g^{ij}\rho_i\rho_i=g^{11}\rho_1^2 \neq 0$ on $\{x_1=0\} \cap B^{+}$. Then the second boundary condition in \eqref{ap2} implies that $u_1=0$ on $\{x_1=0\} \cap B^{+}$. Then we can extend $u$ to the unit ball $B=\{\sum x_i^2 \leq 1\}$ as
\begin{align*}
\left\{ 
\begin{aligned}
\tilde{u}&=u,\, && B^{+}\\ 
\tilde{u}&=0,\, &&B \setminus B^{+} 
\end{aligned}
\right.
\end{align*} such that $\tilde{u} \in W^{1,p}(B)$.

 Since $g$ is $C^k$ on $B^+$, we can extend $g$ to a $C^k$ Riemannian metric on the unit ball $B$. Then we get 
\begin{align*}
|g|^{-\frac{1}{2}}\partial_{i}\left(|g|^{\frac{1}{2}} g^{ij}\tilde{u}_j\right)&=0\end{align*} on $B$ and $\tilde{u}=0$ on the open set $B\setminus B^{+}$. By the weak unique continuation principle of second order elliptic operators, $\tilde{u}$ must vanish in $B$. It implies that $u$ vanishes in $\Omega$. As a harmonic function on $M$, $u$ vanishing in the open set $\Omega \cap M$ implies that $u$ vanishes everywhere in $M$.\\
\end{proof}
\begin{remark}
See Theorem 1.1 in \cite{GL} for the weak unique continuation principle of second order elliptic operators. (Or see Section 1 in \cite{U} where $C^3$ coefficients are assumed.) \end{remark}

\section{Main Results}
 In this section, we will establish generic properties of Steklov eigenfunctions using transversality theorems. 
 
 Let us recall some basic definitions at first.
\begin{defi}\label{rv}Consider a $C^1$ map $F: X\to Y$ between two manifolds. Say $x\in X$ is a regular point of $F$ if $D_xF : T_x (X) \to T_{F(x)}(Y)$ is onto. Say $y \in Y$ is a {\it regular value} if every point $x \in F^{-1}(y)$ is a regular point. 
\end{defi}
\begin{defi}
Say a subset of a topological space is {\it residual} if it is the countable intersection of open dense sets. 
\end{defi} 
\begin{defi}
Say a map $F: X \to Y$ is transverse to a submanifold $Y^{\prime} \subset Y$, if for all $x \in X$ with $F(x) \in Y^{\prime}$, \[
\Ima (D_xF)+T_{F(x)}(Y^{\prime})=T_{F(x)}Y.\]
\end{defi}
 We are going to use transversality theorems from \cite{U} which are stated as follows.
\begin{thm}\label{T1}
Let $\phi: H \times B \to E$ be a $C^k$ map, $H, B$ and $E$ be Banach manifolds with $H$ and $E$ separable. If $0$ is a regular value of $\phi$ and $\phi_b=\phi(, b)$ is a Fredholm map of index less than $k$, then the set $\{b \in B, 0 \, \text{is a regular value of $\phi_b$}\}$ is residual in $B$.
\end{thm}

\begin{thm}\label{T2}
Let $Q, B ,X, Y$ and $Y^{\prime}$ be separable Banach manifolds, $Y^{\prime} \subset Y$ and $Y, Y^{\prime}$ and $X$ finite dimensional. Let $\pi: Q \to B$ be a $C^k$ Fredhlom map of index $0$. Then if $f: Q \times X \to Y$ is a $C^k$ map for $k>\max(1, \rm{dim} X+\rm{dim}Y^{\prime}-\rm{dim} Y)$, and if $f$ is transverse to $Y^{\prime}$, then the set $\{b \in B:  f_b=f|_{\pi^{-1}_{b}} \, \text{is transverse to $Y^{\prime}$ }\}$ is residual in $B$.
\end{thm}

\subsection{Simple \S Eigenvalues}  
 
We will use Theorem \ref{T1} to show that nonzero \S eigenvalues are simple for generic metrics. Precisely, we will prove:

\begin{thm}\label{Gen1}
Let $M^n$ be a smooth compact manifold with smooth boundary $\partial M^n$. Given an integer $k \geq 2$, there exists a residual set of metrics  $g \in \Gk$ for which all nonzero \S eigenvalues are simple.
\end{thm}

Consider $f\in \HM$ with $u \in \wm$ as its harmonic extension. When $f$ is an eigenfunction corresponding to the non-zero \S eigenvalue $\lambda$, it is easy to see that $f$ is non-constant by \eqref{eq1} and \eqref{eq1.1}. Moreover, by Green's formula, \[\int_{\partial M^n}f dA=\frac{1}{\lambda}\int_{\partial M^n}u_n \cdot 1 dA=\frac{1}{\lambda}\int_{M^n} \Delta u \,dV=0.\]

So we introduce the following space in order to apply Theorem \ref{T1}:
 \begin{align*}
 S_k^p&=\left\{f \in \HM: \int_{\partial M}f dA=0\right\}. 
 \end{align*}
And define the map 
\begin{align*}
\phi:  S^p_k \times \R \times \Gk &\rightarrow  \hm\\
(f, \lambda, g)&\rightarrow (\Lambda-\lambda)f.
\end{align*}
When $(f, \lambda, g) \in \phi^{-1}(0)$, it follows that $\lambda$ is a nonzero \S eigenvalue with respect to metric $g$ and $f$ is one corresponding \S eigenfunction which is non-constant.

Here the tangent space can be characterized as
\begin{align*}
&TS^p_k \times \R \times T\Gk|_{(f, \lambda, g)}= \bigg\{
(\nu, s, h): \nu \in \hm, \\
&\int_{\partial  M^n}\nu  dA=0, s \in \R , h=Dg\bigg\}.
\end{align*} 

The differential of $\phi$ is given as follows: 
\begin{align}\label{DP}
D \phi_{(f, \lambda, g)}(\nu, s, h)= (\Lambda -\lambda) \nu-s f+(D_g \Lambda)(h) f
\end{align} with $(\nu, s, h) \in TS^p_k \times \R \times T\Gk|_{(f, \lambda, g)}$.


For convenience, we will denote the differential of $\phi$ in the direction of $S^p_k \times \mathbb{R}$ by $D_1$ and the differential in the direction $\Gk$ by $D_g$. That is,  
\begin{align*}
D\phi=&(D_1 \phi, D_g \phi)\,\, \text{at $(f, \lambda, g)$}\\
\text{with}&\left\{
\begin{aligned}
&D_1 \phi(\nu, s, h)=(\Lambda -\lambda) \nu-s f\\
&D_g \phi(\nu, s, h)=(D_g \Lambda)(h) f.
\end{aligned}
\right.
\end{align*}

We will apply Theorem \ref{T1} to the map $\phi$ introduced above. For this purpose, we need to show that $0$ is a regular value of $\phi$. At first, let us prove the following density lemma about the image of $D_{g} \phi$.

\begin{lemma}\label{density}
At any point $(f, \lambda, g ) \in \phi^{-1}(0)$, there is 
\begin{align*}
( \rm{Im}\,D_g\phi)^{\bot}&\subseteq\{\psi \in \hm \,|\,\rm{supp}(\psi) \subseteq f^{-1}(0) \}.
\end{align*}
In other words, 
\begin{align*}
\{\psi  \in \hm\,|\, \rm{supp}(\psi )\cap f^{-1}(0)=\emptyset\}&\subseteq \rm{Im}(D_g \phi).
\end{align*}Moreover, for any $\psi \in \hm$ with $\rm{supp}(\psi)\cap f^{-1}(0)=\emptyset$, there exists $\sigma \in \Ck$ such that $(D_g \L)(\sigma g)f=\psi$.
\end{lemma}
\begin{proof}
 Fix $(f, g, \lambda)\in \phi^{-1}(0)$. We know that $\lambda \neq 0$ by the definition of $\phi$. Let $u$ be the harmonic extension of $f$ and $v=D_g u$. Take $\psi \in ( \Ima D_g\phi|_{(f, g, \lambda)})^{\bot}$ and denote its harmonic extension by $\tilde{u}$. 
 
 Since $D_g\phi(\nu, s, h)=(D_g \Lambda)(h)f$ at $(f, g, \lambda)$, it follows that for any $h \in T\Gk$:
\begin{align}\label{or}
0&=\int_{\partial  M^n}\psi (D_g \Lambda)(h)f dA.
\end{align}  

When $h=\sigma g$ as a conformal variation for arbitrary $\sigma \in \Ck$, apply Lemma \ref{VL} to \refeq{or}. It follows:
\begin{align*}
0&=\int_{\partial  M^n} \psi (D_g \Lambda)(\sigma g)f \,dA\\
&=-\frac{1}{2}\int_{\partial  M^n}\sigma \psi u_n \,dA+\int_{\partial  M^n}\psi v_n \,dA.
\end{align*}
Then apply Green's formula to the second integral:
\begin{align*}
0&=-\frac{1}{2}\int_{\partial  M^n}\sigma \psi u_n \,dA+\int_{ M^n}(\tilde{u} \Delta v +\langle \nabla \tilde{u}, \nabla v \rangle) \,dV\\
&=-\frac{1}{2}\int_{\partial  M^n}\sigma \psi u_n \,dA+\int_{ M^n}(\tilde{u} \Delta v -v\Delta \tilde{u}) \,dV+\int_{\partial  M^n} v\tilde{u}_n \,dA\\
\end{align*}
Since $\tilde{u}$ is the harmonic extension of $\psi$, there is $\Delta \tilde{u}=0$ on $M^n$. Since $v$ is the variation of $u$, Lemma \ref{VH} implies that $v|_{\partial  M^n}=0$. Thus the equation above becomes :
\begin{align}\label{eq311}
0&=-\frac{1}{2}\int_{\partial  M^n}\sigma \psi u_n \,dA+\int_{ M^n}\tilde{u} \Delta v \,dV.
\end{align}

By Lemma \ref{VH}, we know
\begin{align*}
\Delta v&=(1-\frac{n}{2})\langle \nabla \sigma, \nabla u\rangle.
\end{align*} Apply this to the last integral in \eqref{eq311}. Then we get
\begin{align*}\
0
&=-\frac{1}{2}\int_{\partial  M^n}\sigma \psi u_n\, dA+(1-\frac{n}{2})\int_{ M^n}\tilde{u} \langle \nabla \sigma, \nabla u\rangle \,dV.
\end{align*} 
Apply the Green's formula to the second integral. It follows  
\begin{align}\label{eq312}
0&=-\frac{1}{2}\int_{\partial  M^n}\sigma\psi u_n \,dA+(1-\frac{n}{2})\left(-\int_{ M^n}(\sigma \tilde{u}\Delta u +\sigma \langle \nabla \tilde{u}, \nabla u\rangle)\,dV+\int_{\partial  M^n}\sigma \psi u_n \,dA\right).
\end{align}

Since 
 \begin{align*}
 \left\{
\begin{aligned}
\Delta u&=0, \, &&M^n,\\
u_n|_{\partial  M^n}&= \lambda f,\, &&\partial  M^n,
\end{aligned}%
\right.
\end{align*} 
the equation \eqref{eq312}  can be simplified as
\begin{align}\label{eq20}
0&=-\frac{\lambda(n-1)}{2}\int_{\partial  M^n}\sigma\psi f dA-(1-\frac{n}{2})\int_{ M^n}\sigma \langle \nabla \tilde{u}, \nabla u\rangle dV.
\end{align}

 Since $\sigma \in \Ck$ is arbitrary and $\lambda \neq 0$, the equation \eqref{eq20} implies that $\psi f=0$ on $\partial  M^n$. This means :
\begin{align*} 
\Ima^{\bot}(D_g\phi)&\subseteq \{(D_g \L)(\sigma g)f\,|\, \sigma \in \Ck\}^{\bot}\\
&\subseteq\{\psi  |\, \psi f=0\, \text{ on $\partial M^n$}\}\cap \hm\\
&=\{\psi |\, {\rm supp}\, \psi \subseteq f^{-1}(0)\} \cap \hm.
\end{align*}
Therefore we get 
\begin{align*}
\Ima (D_g \phi)&\supseteq \{(D_g \L)(\sigma g)f\,|\, \sigma \in \Ck\}\\
&\supseteq \{\psi\,|\, {\rm supp}(\psi)\cap f^{-1}(0)=\emptyset\}\cap \hm.
\end{align*}
The second relation implies that for any $\psi \in \hm$ with $\rm{supp}(\psi)\cap f^{-1}(0)=\emptyset$, there exists $\sigma \in \Ck$ such that $(D_g \L)(\sigma g)f=\psi$.
\end{proof}

\begin{remark}\label{dense}By the weak unique continuation principle (Theorem \ref{wucp}), $f^{-1}(0)$ is closed with empty interior. So 
$ \{\psi\,| \,\rm{supp}(\psi)\cap f^{-1}(0)=\emptyset\}\cap \hm
$ is dense in $\hm$.
\end{remark}

\begin{prop}\label{reg}
$0$ is a regular value of $\phi$.
\end{prop}
\begin{proof}
By Definition \ref{rv}, we need to prove  $ \Ima D\phi|_{(f, \lambda, g)} = \hm$ at any point $(f, \lambda, g) \in \phi^{-1}(0)$.

We first show that $ \Ima D\phi|_{(f, \lambda, g)}$ is dense in $\hm$. 
Take $\psi\in(\Ima D\phi_{(f, \lambda, g)})^{\bot}$. Then for any $(v,s,h) \in TS^p_k \times \R \times T\Gk|_{(f, \lambda, g)}$, there is 
\begin{align}\label{orth}
0&=\int_{\partial M}\psi\left((\Lambda -\lambda) \nu-s f+(D_g \Lambda)(h) f\right) dA.
\end{align}

When $\nu=0$ and $s=0$ in \eqref{orth}, it follows that for any $h\in T\Gk$
\begin{align*}
0&=\int_{\partial M}\psi\left((D_g \Lambda)(h) f\right) dA.
\end{align*} So $\psi \in (\Ima D_g(\Lambda)f)^{\bot}$. Then Lemma \ref{density} implies that $\rm{supp} \,\psi \subseteq f^{-1}(0)$. By Theorem \ref{wucp}, $f^{-1}(0)$ is closed with empty interior. Therefore $\psi=0$ on some open set of $\partial M^n$.

When $s=0$ and $h=0$ in \eqref{orth}, it follows that for any $\nu\in \HM$ satisfying $\int_{\partial M^n}\nu  dA=0$, there is
\begin{align*}
0&=\int_{\partial M^n}\psi\left((\Lambda -\lambda) \nu\right) dA.
\end{align*} That is, $\psi \in \left(\Ima (\Lambda-\lambda)\right)^{\bot}\cap \hm$. Since $\Lambda-\lambda$ is self-adjoint and elliptic, there is $\hm=\Ima (\Lambda-\lambda)\oplus \ker (\Lambda-\lambda)$. Therefore, we have $\psi \in \ker (\Lambda-\lambda)$. At the same time, the previous argument implies that $\psi=0$ on some open set of $\partial M^n$. By Theorem \ref{wucp}, $\psi$ mush vanish everywhere in $\partial M^n$ as in the kernel of $\Lambda -\lambda$. This means $(\Ima D\phi|_{(f, \lambda, g)})^{\bot}=\{0\}$, i.e., $ \Ima D\phi|_{(f, \lambda, g)}$ is dense in $\hm$.


 Since $\L-\lambda$ is Fredholm, it follows that $\Ima (\Lambda-\lambda)$ is closed in $\hm$ with the finite co-dimension denoted by $l$. By the density of the image of $D\phi|_{(f, \lambda, g)}$, we can select a $l-$dimensional subspace $V \subset T\Gk|_{(f, \lambda, g)}$ such that 
\[
\hm= \Ima (\Lambda-\lambda) \oplus D\phi|_{(f, \lambda, g)}(V) \subseteq \Ima  D\phi|_{(f, \lambda, g)}.
\]Hence $0$ is a regular value of $\phi$.
\end{proof}

Now we are ready to prove Theorem \ref{Gen1}.

\begin{proof}[Proof of Theorem \ref{Gen1}] We apply Theorem \ref{T1} with $H=S^p_k \times \R, B=\Gk, E=\hm$ here. 

By Theorem \ref{reg}, we know that $0$ is a regular value of $\phi$. We need to show that $\phi_g$ is a Fredholm map of index $0$. As a self-adjoint elliptic operator, $\L-\lambda: \HM \rightarrow \hm$ has the Fredhlom index $0$ for a fixed $\lambda$. The restriction on $S^p_k$ will reduce the index of $\Lambda-\lambda$ by one since it has the co-dimension one in $\hm$. So the index of $\phi_g=\L-\lambda$ will be zero when $\lambda$ is allowed to vary. Then Theorem \ref{T1} implies that the set of $g \in \Gk $ such that $\phi_g$ with $0$ as a regular value is residual.

The conclusion follows since that $\L$ with respect to $g$ has all nonzero eigenvalues simple if and only if $\phi_g$ has $0$ as a regular value. To show this, it is enough to prove that:\\
{\it For a fixed $g$, all non-zero eigenvalue of $\Lambda$ are simple if and only if $D\phi_g$ is onto at every $(f, \lambda) \in \phi^{-1}_g(0)$.}\\
We notice that $(f, \lambda) \in \phi^{-1}_g(0)$ if and only if $f$ is one eigenfunction of $\Lambda$ corresponding to $\lambda \neq 0$. 
At $(f, \lambda) \in \phi^{-1}_g(0)$, by \eqref{DP}, there is 
\begin{align*}
D\phi_g : TS^p_k \times \R &\rightarrow \hm\\
(v, s)&\rightarrow (\Lambda -\lambda) v-sf.
\end{align*}
 Then it follows \[\Ima D\phi_g=\Ima (\Lambda -\lambda) \oplus \rm{span}(f).\] On the other hand, since $\Lambda-\lambda$ is elliptic and self adjoint, we have 
 \[\hm=\Ima (\Lambda -\lambda)\oplus \ker (\Lambda -\lambda).\] Therefore, at $(f, \lambda) \in \phi^{-1}_g(0)$, $D\phi_g$ is onto $\hm$ if and only if $\ker (\Lambda -\lambda)=\rm{span}(f)$, that is, $\lambda$ is simple.

%

\end{proof}

\subsection{Nodal sets and critical points of \S eigenfunctions}

We will apply Theorem \ref{T2} to study nodal sets and critical points of non-constant \S eigenfunctions.  We will work on the Banach manifold of all non-constant eigenfunctions $Q=\phi^{-1}(0)$ which is a smooth manifold by Proposition \ref{reg}. Let $\pi: Q \to \Gk$ be the restriction of the projection $S^p_k\times \R \times \Gk$ to $\Gk$ on the parameter space. Then $\pi$ is a Fredholm map of index $0$ from the dimension counting.

Let us define the evaluation maps as:
 \begin{align*}
 \alpha: Q \times \partial  M^n &\rightarrow \mathbb{R}, \\
 (f, \lambda, g, x) &\rightarrow f(x),
 \end{align*}
 and 
 \begin{align*}
  \beta: Q \times \partial  M^n &\rightarrow T^{\ast}(\partial  M^n),\\
  (f, \lambda, g, x) &\rightarrow df(x).
\end{align*}Here the notation $d$ is the exterior derivative with respect to $\partial M^n$ instead of $M^n$.

To apply Theorem \ref{T2} to the map $\alpha$, we need to show it has $0$ as a regular value at first.

\begin{prop}\label{regA}
$0$ is a regular value of $\alpha$.
\end{prop}
\begin{proof} For $x \in \partial M^n$, let $\alpha_x =\alpha(, x): Q \rightarrow \R$. To show that $0$ is a regular value of $\alpha$, it is enough to show for every $x \in \partial  M^n$, 
\begin{align*}
D\alpha_x: T_{(f, \lambda, g)} Q &\rightarrow \R, \\
(v, s, h) &\rightarrow v(x)
\end{align*}  is onto when $\alpha_x (f, \lambda, g)=f(x)=0$. 

As a linear map, $D\alpha_x$ is either onto $\mathbb{R}$ or $D\alpha_x (v, s, h)=v(x)=0$ for all $(v, s, h) \in T_{(f, \lambda, g)} Q$. Assume the second case holds for some $x \in \partial M$. We will construct some special $(\nu, 0, h) \in T_{(f, \lambda, g)} Q$ to get a contradiction. 

Notice the constraint equation for $(\nu, s, h) \in T_{(f, \lambda, g)} Q$ is :
\begin{align}\label{CE}
\left\{
\begin{aligned}
&(\Lambda-\lambda)\nu-sf+D_g(\Lambda)(h)f=0,\\
&\int_{\partial  M^n} \nu \,dA=0.
\end{aligned}
\right.
\end{align}

Choose $w \in \ker^{\bot}(\Lambda -\lambda)$ with $\rm{supp}(w)\cap f^{-1}(0)=\emptyset$ and $\int_{\partial M}w dA=0$. According to Lemma \ref{density}, there exists $\sigma_w\in \Ck$ such that $h_w=\sigma_w g$ satisfies
\[D_g(\Lambda)(h_w)f= -w.\]

Since $\Lambda$ is self-adjoint, there exists sequence of eigenfunctions $\{\psi_n\}$ corresponding to eigenvalues $\lambda_n$ which forms an orthonormal basis of $L^2(\partial M)$.  As discussed in Section 3.1, when the eigenvalue $\lambda_n=0$, the corresponding eigenfunction is constant. When the eigenvalue $\lambda_n>0$, the corresponding eigenfunction $\psi_n$  is non-constant and satisfies $\int_{\partial M} \psi_n dA=0$ (i.e. $\psi_n$ is orthogonal to $1$). 

 Since $\Lambda$ has a compact resolvent (\cite{AEKS}), we can consider the resolvent operator 
\[R_{\lambda} w=\underset{ \lambda_n\neq \lambda}{\sum} \frac{\langle w, \,\phi_n\rangle_{\partial M}}{\lambda_n-\lambda}\psi_n\] which satisfies 
\[
(\Lambda-\lambda)R_{\lambda} w=w.
\] Here $\langle w, \,\phi_n\rangle_{\partial M}=\int_{\partial M}w\psi_n dA$. In addition, $\int_{\partial M}w dA=0$ implies that 
\[
R_{\lambda} w=\underset{ \lambda_n \neq \lambda, \lambda_n>0}{\sum} \frac{\langle w, \,\phi_n\rangle_{\partial M}}{\lambda_n-\lambda}\psi_n.\]

Then $(R_{\lambda} w, 0, h_w) \in T_{(f, \lambda, g)}Q$ according to \eqref{CE}. By the assumption, $R_{\lambda} w(x)=0$ for all such admissible $w$, i.e.,
\begin{align}\label{rw}
\underset{\lambda_n \neq \lambda, \lambda_n>0}{\sum} \frac{\langle w, \psi_n\rangle_{\partial M}}{\lambda_n-\lambda}\psi_n(x)=0.
\end{align}

According to Remark \ref{dense}, the set $\{w\,|\,\rm{supp}(w)\cap f^{-1}(0)=\emptyset \}$ is dense in $\HM$. So the set of these admissible $w$ is dense in $ \ker^{\bot}(\Lambda-\lambda) \cap \{\psi\in \HM)|\int_{\partial M}\psi dA=0\}$. Then \eqref{rw} implies that
\begin{align}\label{321}
\psi_n(x)=0 \, \text{when $\lambda_n \neq \lambda, \lambda_n>0$.}
\end{align}

When $\psi_n$ is the eigenfunction of $\lambda$, we get $(\psi_n, 0, 0)\in T_{(f, \lambda, g)}Q$ according to \eqref{CE}. Therefore $\psi_n(x)=0$ by the assumption when $\lambda_n=\lambda$. This fact together $\eqref{321}$ imply that all non-constant eigenfunctions of $\Lambda$ with respect to $g$ vanish at $x$ which is absurd. Thus $D\alpha_x$ is onto $\R$.
\end{proof}

\begin{thm}\label{Gen2}
Let $M^n$ be a compact manifold with smooth boundary $\partial  M^n$ and $n\geq 2$. Given an integer $k> n-1$, there exists a residual set of metrics $g\in \Gk$ such that for all $g$ in this set, non-constant \S eigenfunctions have zero as a regular value in $\partial M^n$.
\end{thm}
\begin{proof} 

We apply Theorem \ref{T2} here with replacing $f$ by $\alpha$, $B$ by $\Gk$, $X$ by $\partial M^n$, $Y$ by $\R$ and $Y^{\prime}$ by $\{0\}$.

By the Sobolev imbedding theorem, there is $\HM \subset C^{k-1}(\partial  M^n)$ for any $p>n$. So $\alpha$ is a $C^{k-1}$ map for any $p>n$ which is sufficient smooth to apply Theorem \ref{T2}, with $\rm{dim} X+\rm{dim} Y^{\prime}-\rm{dim} Y=n-2$ and the assumption $k>n-1$. Moreover, $\alpha$ has $0$ as a regular value by Proposition \ref{regA}. Therefore Theorem \ref{T2} implies that there exits a residual set of $g \in \Gk$ such that $\alpha_g=\alpha|_{\pi^{-1}(g)}$ has $0$ as a regular value. 

We know that $\alpha_g(f, \lambda, x)=0$ if and only if $f$ is one eigenfunction of $\Lambda$ corresponding to $\lambda \neq 0$ with respect to $g$ and $f(x)=0$. And $d\alpha_g=df|_{x}$ at $(f, \lambda, x)$. Therefore $\alpha_g$ has $0$ as a regular value implies that that every non-constant \S eigenfunction $f$ of $\Lambda$ with respect to $g$ has $0$ as a regular value. The conclusion follows.

\end{proof}

Next we will apply Theorem \ref{T2} to the map $\beta$. We need to prove the following proposition at first.
\begin{prop}\label{regB}
$\beta$ is transverse to the zero section of $T^{\ast}(\partial  M^n)$.
\end{prop}
\begin{proof}
For $x\in \partial M^n$, let $\beta_x=\beta(, x): Q \rightarrow T_x^{\ast}(\partial  M^n)$. To show that $\beta$ is transverse to the zero section of $T^{\ast}(\partial  M^n)$, it is enough to show that for every $x\in \partial M^n$,
\begin{align*}
D\beta_x : T_{(f, \lambda, g)} Q  &\rightarrow T_x^{\ast}(\partial M^n)\\
 (\nu, s, h) &\rightarrow d\nu(x)
\end{align*} is onto when $\beta_x(f, \lambda, g)=df(x)=0$.

Assume this is not true for some $(f,\lambda, g)$ and $x$ with $df(x)=0$. Since $\rm{Im} D\,\beta_x$ is a linear subspace of $T_x^{\ast}(\partial M^n)$, there will exist a nonzero $\xi \in T^{\ast}(\partial M^n)$ such that $\xi \cdot d\nu(x) =0$ for all $d\nu(x) \in \rm{Im}\, D\beta_x$ by the assumption. Here $\xi \cdot d\nu(x)$ denotes the inner product in $T^{\ast}(\partial M^n)$. 

Consider the special $(\nu, 0, h)$ constructed in the proof of Proposition \ref{regA}. That is, for any $w \in \rm{Ker}^{\bot}(\Lambda -\lambda)$ with $\rm{supp}(w)\cap f^{-1}(0)=\emptyset$ and $\int_{\partial M} w dA=0$, there is 
\[\nu=R_{\lambda} w=\underset{\lambda_n \neq \lambda, \lambda_n>0}{\sum} \frac{\langle w, \,\psi_n\rangle_{\partial M}}{\lambda_n-\lambda}\psi_n\] 
and 
\[h_w=\sigma_w g\, \text{for some $\sigma_w\in \Ck$}\]  such that $(R_{\lambda} w,0,h_{w}) \in T_{(f, \lambda, g)}Q$. Then for all such admissible $w$, it follows that
\[
\xi \cdot d R_{\lambda}w(x)=0.\]
That is,
\[ \underset{\lambda_n \neq \lambda, \lambda_n>0}{\sum} \frac{\langle w, \psi_n\rangle_{\partial M}}{\lambda_n-\lambda}\xi \cdot d\psi_n(x)=0.\] Then the density of the set of the admissible $w$ in $\rm{Ker}^{\bot}(\Lambda -\lambda)\cap \{\int_{\partial M} \psi dA=0\}$ implies that 
\begin{align}\label{eigen}
&\xi\cdot d\psi_n(x)=0\, \, \text{when $\lambda_n \neq \lambda, \lambda_n>0$.}
\end{align}

When $\psi_n$ is the eigenfunction of $\lambda$, we get $(\psi_n, 0, 0)\in T_{(f, \lambda, g)}Q$ according to \eqref{CE}. Then $\xi\cdot d\psi_n(x)=0$ by the assumption when $\lambda_n=\lambda$. In addition, when $\psi_n$ is the eigenfunction of the eigenvalue $\lambda_n=0$, $\psi_n$ is constant and then $\xi\cdot d\psi_n(x)=0$. These facts together with $\eqref{eigen}$ imply that $\xi \cdot d\zeta(x) =0$ for any function $\zeta$ on $\partial M^n$. For any $\eta \in T_x^{\ast}(\partial M^n)$, we can find some function $\zeta$ on $\partial M^n$ such that $\eta=d \zeta$ at the point $x$. It follows that $\xi \cdot \eta(x) =0$ for all $\eta \in T_x^{\ast}(\partial M^n)$. Then $\xi =0$ which is a contradiction. Therefore $D \beta_x$ is onto.

\end{proof}

\begin{thm}\label{Gen3}
Let $M^n$ be a compact manifold with smooth boundary $\partial M^n$. Given an integer $k \geq 2$, there exists a residual set of metrics $g\in \Gk$ such that for any $g$ in this set, all non-constant \S eigenfunctions are Morse functions.
\end{thm}

\begin{proof}
We apply Theorem \ref{T2} here with replacing $f$ by $\beta$, $B$ by $\Gk$, $X$ by $\partial M$, $Y$ by $T^{\ast}(\partial M^n)$ and $Y^{\prime}$ by the zero section of $T^{\ast}(\partial M^n)$.

Since $f \in \HM \subset C^{k-1}(\partial M^n)$ for all $p\geq n$, we have $\beta$ as a $C^{k-1}$ map for all $p\geq n$. So $\beta$ is sufficiently smooth to apply Theorem \ref{T2}, with $\rm{dim} X+\rm{dim} Y^{\prime}-\rm{dim} Y=0$ and the assumption $k\geq 2$. By Proposition \ref{regB}, $\beta$ is transverse to $Y^{\prime}$, the zero section of $T^{\ast}(\partial M^n)$. Thus Theorem \ref{T2} implies that there exits a residual subset of $g \in \Gk$ such that $\beta_g=\beta|_{\pi^{-1}(g)}$ is transverse to the zero section of $T^{\ast}(\partial M^n)$.

We know that $\beta_g(f, \lambda, x)=0$ if and only if $f$ is one eigenfunction of $\Lambda$ corresponding to $\lambda \neq 0$ with respect to $g$ and $df(x)=0$. Therefore $\beta_g$ is transverse to the zero section of $T^{\ast}(\partial M^n)$ will imply that the map $df :\partial M^n \to T^{\ast}(\partial M^n)$ is transverse to the zero section of $T^{\ast}(\partial M^n)$ for any non-constant eigenfunction $f$ with respect to $g$. This means that in a local coordinates $U \subset \partial M^n$, the map 
\[
(f_{x_1}, \cdots, f_{x_{n-1}}): U \to \R^{n-1}
\] has $0$ as a regular value, which is equivalent to the Hessian of $f$ at critical points being non-degenerate. Thus the conclusion follows.

\end{proof}

In the end, Theorem \ref{Gen1}, \ref{Gen2} and \ref{Gen3} together imply the Main Theorem.

\appendix

\end{document}